\documentclass{article}

\usepackage{arxiv}

\usepackage[utf8]{inputenc} % allow utf-8 input
\usepackage[T1]{fontenc}    % use 8-bit T1 fonts
\usepackage{hyperref}       % hyperlinks
\usepackage{url}            % simple URL typesetting
\usepackage{booktabs}       % professional-quality tables
\usepackage{amsfonts}       % blackboard math symbols
\usepackage{nicefrac}       % compact symbols for 1/2, etc.
\usepackage{microtype}      % microtypography
\usepackage{lipsum}		% Can be removed after putting your text content
\usepackage{graphicx}
\usepackage{natbib}
\usepackage{doi}

\title{Energy Efficient Manufacturing Scheduling: A Systematic Literature Review}

\author{ {\includegraphics[scale=0.06]{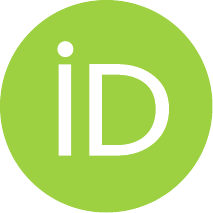}\hspace{1mm}Ahmed Missaoui}\thanks{Use footnote for providing further
		information about author (webpage, alternative
		address)---\emph{not} for acknowledging funding agencies.} \\
	Department of Computer Science\\
	University College Cork\\
	Cork, T12 XF62, Ireland \\
	\texttt{amissaoui@ucc.ie} \\
\And
	           {\includegraphics[scale=0.06]                       
        {orcid.pdf}\hspace{1mm}Cemalettin Ozturk} \\
	Process, Energy and Transport Engineering\\
	  Munster Technological University\\
	Bishopstown, Cork, T12 P928, Ireland\\
	\texttt{cemalettin.ozturk@mtu.ie} \\
 \And
       {\includegraphics[scale=0.06]{orcid.pdf}\hspace{1mm}Barry O’Sullivan} \\
		Department of Computer Science\\
	University College Cork\\
	Cork, T12 XF62, Ireland \\
	\texttt{b.osullivan@cs.ucc.ie} \\
 \And
        {\includegraphics[scale=0.06]{orcid.pdf}\hspace{1mm}Michele Garraffa} \\
		Department of Computer Science\\
	University College Cork\\
	Cork, T12 XF62, Ireland \\
	\texttt{michele.garraffa@gmail.com} \\
}
\begin{document}
\maketitle

\begin{abstract}
The social context in relation to energy policies, energy supply, and sustainability concerns as well as advances in more energy-efficient technologies is driving a need for a change in the manufacturing sector. 
The main purpose of this work is to provide a research framework for energy-efficient scheduling (EES) which is a very active research area with more than 500 papers published in the last 10 years. The reason for this interest is mostly due to the economic and environmental impact of considering energy in production scheduling. In this paper, we present a systematic literature review of recent papers in this area, provide a classification of the problems studied, and present an overview of the main aspects and methodologies considered as well as open research challenges.
\end{abstract}

% keywords can be removed

\section{Introduction}
\label{sec:Introduction}
In 2021, the manufacturing industry accounted for 25.6 \% of the energy consumed within the EU (\cite{Euro2021}). Besides, European countries rely on imports to meet their energy needs, with an average energy dependency rate of 60.7 \%.
This aspect constitutes a major threat to energy security and can be mitigated by reducing overall energy consumption. While there are substantial amounts to reduce the energy consumption of manufacturing equipment and machining, research shows that the majority of the energy waste is observed during machine idle time and/or setup times. Hence, it is obvious that system view and implementing intelligent manufacturing scheduling methods have the potential to help to reduce energy consumption in the manufacturing industry.  
One direction to achieve this goal is exploiting Operations Research (OR) and Artificial Intelligence (AI) techniques to optimize the production schedules in the industry. The basic idea of these methods is to reduce the energy cost of the system while still achieving a good level of productivity.
An extensive amount of research has been recently conducted to propose the use of these techniques in a heterogeneous set of applications in production scheduling, arising in very different countries around the world from different industries.
As an example: methods for energy-efficient production scheduling of tortilla manufacturing are proposed in \cite{yaurima2018hybrid}, while a case study from a Chinese motor company is studied in \cite{Lu2017}. An Iranian extractor hood factory, \cite{Ramezanian2019}, and a Belgian plastic bottle manufacturer, \cite{Gong2017}, are two other examples. 

In this work, we review the recent research in the field of energy-efficient production scheduling, by analyzing and classifying the literature published since 2013.
The key message is that the number of papers published in this area has been sharply increasing, in fact, the value in 2021 is more than triple of the value reached in 2013.
The reason for that is surely linked to economic and political aspects as stated above, but also due to the challenging and rewarding scientific nature of the area. %we also believe that researchers have been attracted by this area because this type of production scheduling problems are very challenging and rewarding from a scientific perspective.
There are various literature review papers on the energy-efficient manufacturing domain 
(give citations to them).
\cite{Renna2021}
\cite{Terbrack2021}
However, our work is contributing to the area by (1) focusing on the latest trends in the last decade (2) providing a holistic view of the characteristics of the problem both in terms of objective functions, problem formulation, and solution method (3) comparing to some of the review papers that are limited to a particular production systems 
( \cite{Para2022} Jobshop,  
\cite{utama2023systematic} hybrid flowshop )
To the best of our knowledge, this is the first systematic literature review that focuses on evaluating the energy efficiency of scheduling strategies across various production systems taking into account all sustainable related issues. 
By considering factors such as system design, workload characteristics, and resource utilization, our literature aims to provide a comprehensive and objective assessment of energy-efficient scheduling strategies in production systems.

The rest of the paper is structured as follows. Section \ref{sec:Review methodology} describes in detail the process of research and the review methodology.
Section \ref{sec:results} discusses general findings and formal results.
Section \ref{sec:classification} presents in detail the classification of problems, describes the different dimensions we use for classifying the literature and provides statistical insights about the most studied settings.
Section \ref{sec:formulation} addresses the different methods used to formulate problems as well as a survey of the most frequently used techniques to solve including both exact and heuristic techniques.
Section \ref{sec:interaction} introduces the relation between objective functions and energy-related aspects. 
Finally, Section \ref{sec:futurework} briefly describes directions for future studies in this area and provides some conclusions. We provide a nomenclature in Table \ref{table:nomenclature} for the rest of the paper.

\begin{table}[htbp]
  \centering
  \begin{tabular}{ | l | p{9cm} |}
    \hline
    \textbf{Abbreviation} & \textbf{Definition} \\
    \hline
    FS & Flowshop \\
    JS & Jobshop \\
    PM & Parallel machines \\
    SM & Single machine \\
    Dis.JS & Distributed Jobshop \\
    Dis. FS & Distributed flowshop \\
    Dis. HFS & Distributed hybrid flowshop \\
    Dis. PM & Distributed parallel machine \\
    HFS & Hybrid flowshop \\
    EES & Energy efficient scheduling \\
    TEC & Total energy consumption \\
    TECost & Total energy cost \\
    Cmax & Makespan \\
    E/T & Earliness tardiness \\
    Car.E & Carbon emission \\
    Pr Cost & Production cost \\
    PBM & Population-based metaheuristics \\
    SSBM & Single solution based metaheuristics \\
    Processing Eng & Processing Energy \\
    Idle Eng & Idle Energy \\
    M.speed & Machine speed \\
    TOU & Time of use \\
    Setup Eng & Setup Energy \\
    Transp Eng & Transportation Energy \\
    T.on/off & Turn on/off \\
    Peak P & Peak Power \\
    Auxiliary Eng & Auxiliary Energy \\
    Maintenance Eng & Maintenance Energy \\
    Renewable Eng & Renewable Energy \\
    Load/unload Eng & Load/unload Energy \\
    Blocking Eng & Blocking Energy \\
    Storage Eng & Storage Energy \\
    Tiered prices & Tiered prices \\
    Preemption Eng & Machine Preemption \\
    \hline
  \end{tabular}
  \caption{Nomenclature}
  \label{table:nomenclature}
\end{table}

\section{Review methodology}
\label{sec:Review methodology}
In this part, we will outline how we conducted the literature review and how we analyzed the recurring themes in the papers we reviewed.
\subsection{Research process}
\begin{figure}[hbtp] 
\centering
\includegraphics[width=1.0\textwidth]{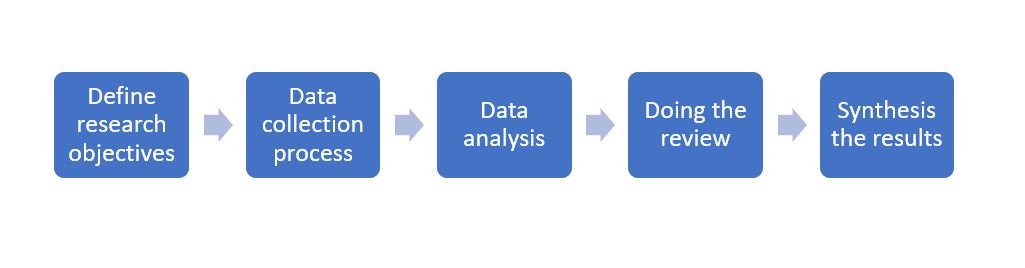} 
\caption{Search steps}
\label{fig:steps}
\end{figure}

In this paper, we follow the systematic review procedures proposed by \cite{andrei2022knowledge} that is summarized in Fig \ref{fig:steps}.

\noindent Reflection of these steps in this literature review is given below: \\
\noindent Step 1: Our research objectives are reviewing the manufacturing scheduling literature in the last decade, providing statistical insights to understand trends both in experimental and applied research, and finally highlighting future research directions.\\  
Step 2: Create search strings and inclusion criteria and by choosing respective abstract and citation databases, carry out the data collection process.

The strategy followed is to look for related publications in the two online abstract and citation databases Scopus and Web of Science. The reason to limit this survey with these two databases are (1) indexing the majority of the publications with high impact in the field (2) allowing fast and precise searches. Search strings used for reviewing the databases are presented in Table \ref{String} which are executed on "article", "book chapter" and "conference paper" types of publications by including "article titles", "abstracts", and "keywords".  

\begin{table}[htbp]
\begin{center}
\begin{tabular}{ | l | p{11cm} |}
\hline
Database & Search strings \\ 
\hline
Scopus  & (TITLE-ABS-KEY("manufacturing" OR "production") AND TITLE-ABS-KEY ("energy" AND "scheduling") AND  TITLE-ABS-KEY ("flowshop" OR "parallel machine" OR "jobshop" OR "hybrid flowshop" OR "flow shop" OR "hybrid flow shop"  OR  "job shop" OR "flow line" OR "assembly line" OR "single")) AND ( LIMIT-TO ( DOCTYPE, "ar" ) OR LIMIT-TO ( DOCTYPE , "cp" ) OR LIMIT-TO ( DOCTYPE , "ch" ) ) AND time period (2013-2023) AND ( LIMIT-TO ( LANGUAGE, "English" )) \\ 
\hline
WoS & "manufacturing" OR "production"(All Fields) AND "energy" AND "scheduling"(All Fields)and"flowshop" OR "parallel machine" OR "jobshop" OR "hybrid flowshop" OR "flow shop" OR "hybrid flow shop" OR "job shop" OR "assembly line" OR "single" OR "flow line"(All Fields)and English (Languages) )  AND  time period (2013-2023) \\ 
\hline
\end{tabular}
\caption{Search Strings}
\label{String}
\end{center}
\end{table}

Figure \ref{fig:search process} describes the search process. After executing the search strings in Table \ref{String} in the first step while Scopus returned 691 documents, WoS provided 757 documents. In the following step, duplicates are removed as the first filtering in the research process which resulted in 931 publications. 
The papers belonging to subject areas outside the scope were excluded after the initial screening of the titles and abstracts, leaving 693 articles in total.\\
Step 3: Conduct the review based on a systematic content analysis following four basic steps; material gathering, assessing the material using an analytical framework, developing classes, and evaluating the material. Developing categories and identifying classes is the antecedent step to reviewing papers and handling inductively where papers were read over and over again. During this process, all used energy aspects, mathematical formulation, shop floor categories, objective functions, and solution methods are established. Finally, all similar characteristics are organized into 15 overarching themes that makes further analysis easier.\\
Step 4: Create a data extraction form that includes all specific criteria before starting the encoding of papers. Using Microsoft Excel the 506 papers are encoded according to 15 technical and non-technical criteria such as year, journal, country, production system, energy-related aspects, and solution method.\\
Step 5: Incorporate and discuss the main outcome of the literature review which will be discussed in detail in Sections\\
\ref{sec:results},\ref{sec:classification},\ref{sec:formulation} and \ref{sec:interaction}.

\begin{figure}[hbtp] 
\centering
\includegraphics[width=1.0\textwidth]{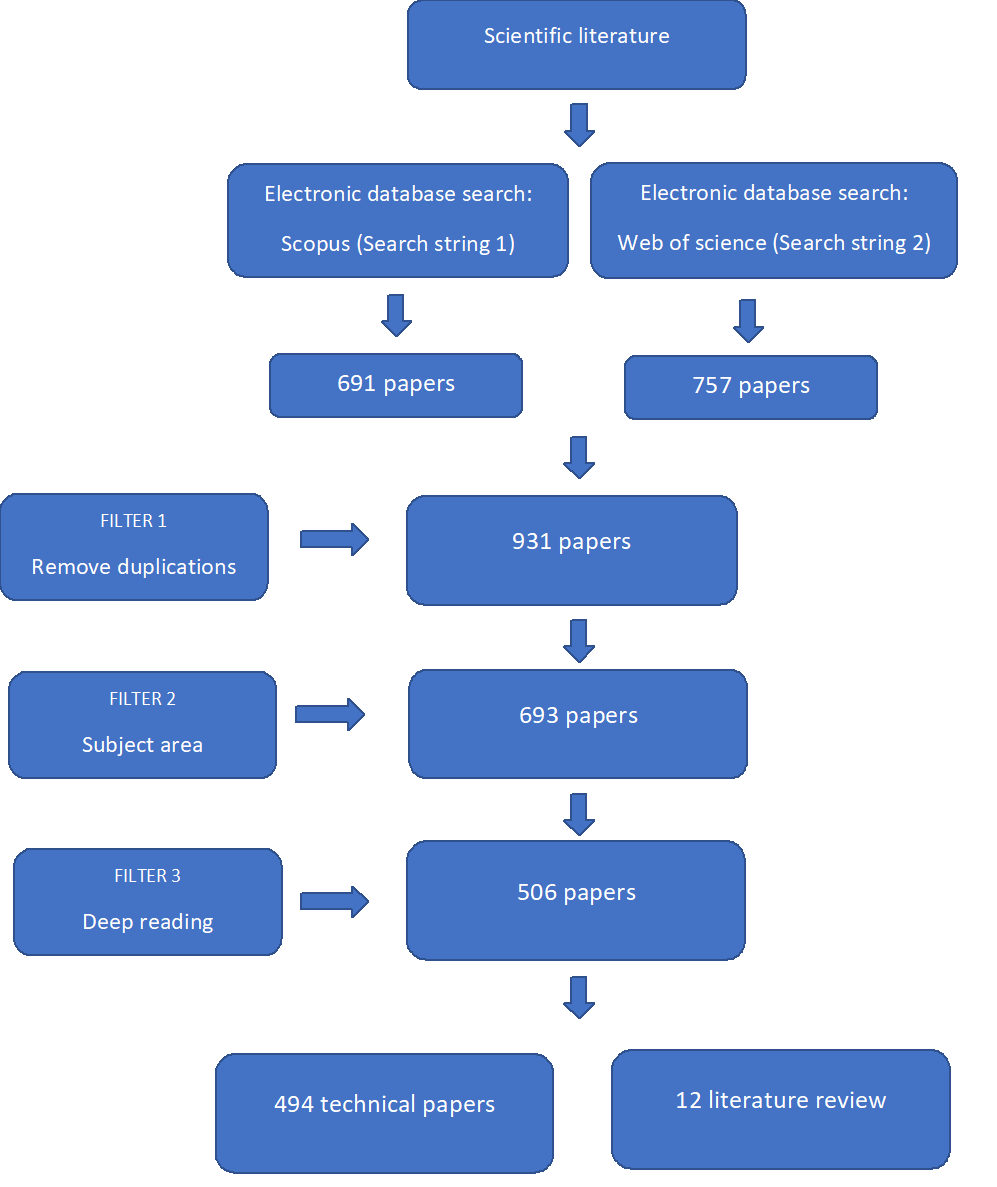}
\caption{Search process}
\label{fig:classes}
\end{figure}

\subsection{Thematic analysis}

\begin{figure}[hbtp] 
\centering
\includegraphics[width=1.0\textwidth]{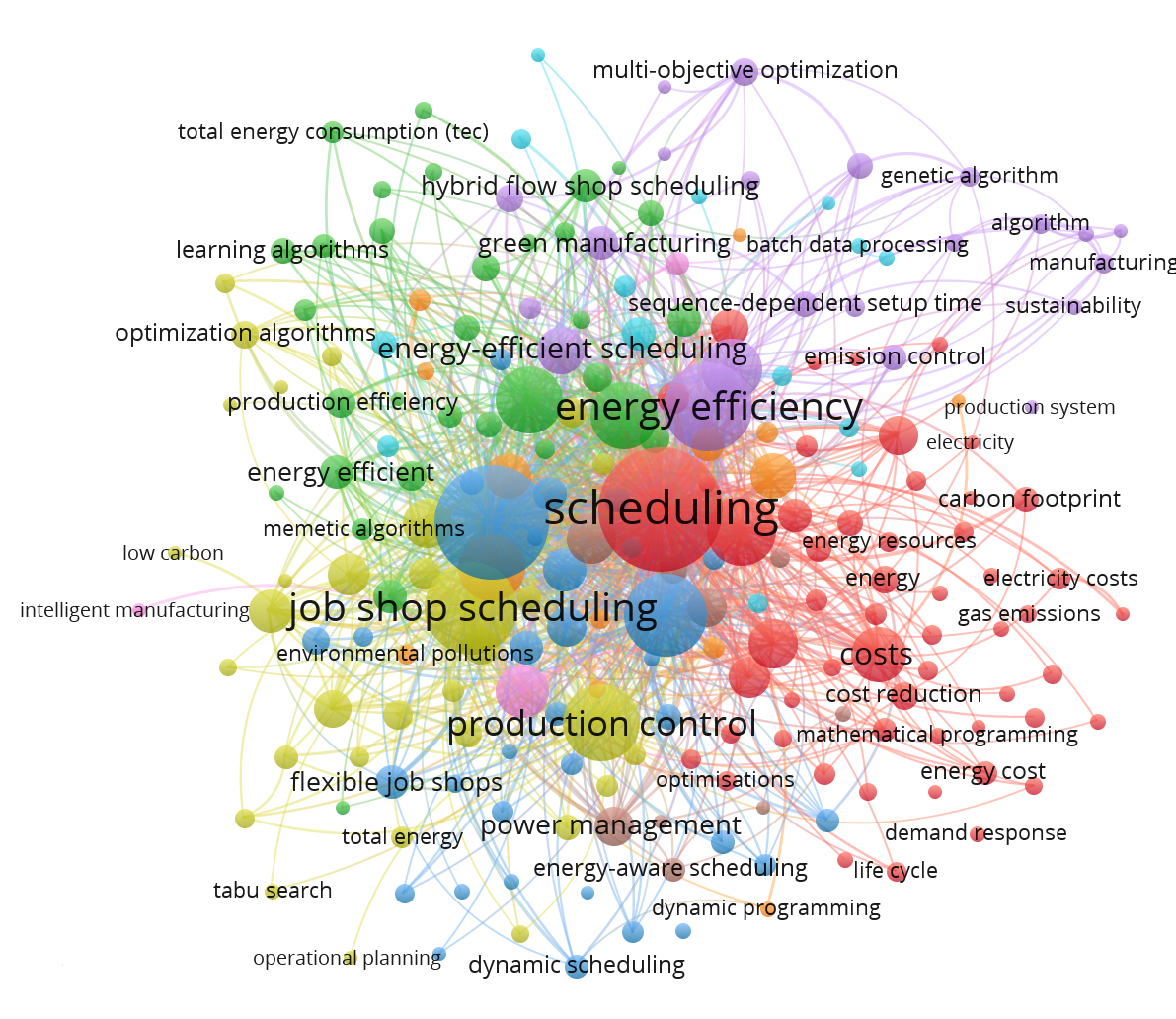}
\caption{Co-occurrence analysis map.}
\label{fig:search process}
\end{figure}

Using the VOSviewer v.1.6.18 tool, a map of co-occurring topics was created based on the keywords of the reviewed publications. The objective in this thematic analysis is to show relationships between any two keywords that exist more than 5 times.
Based on the analysis the most extensively used keyword is "scheduling" which heads a cluster (shown with red circles and edges in Figure \ref{fig:search process}), which groups many objective function considerations related to energy (costs, carbon emission, etc.)   
Another important cluster is shown in blue which is energy utilization which reveals a relationship with some keywords such as "sequence-dependent setup time", "sustainability", and "manufacturing" which justify the systematic view of our review. The thematic analysis also shows that mostly studied manufacturing system is "job shop scheduling" which is directly related to some other keywords relevant to our review such as  "intelligent manufacturing" and "low carbon". 

\section{Results of the systematic literature review}
\label{sec:results}
\subsection{Publication distribution along the years}
\begin{figure}[hbtp] 
\begin{center}

\includegraphics[width=1.0\textwidth]{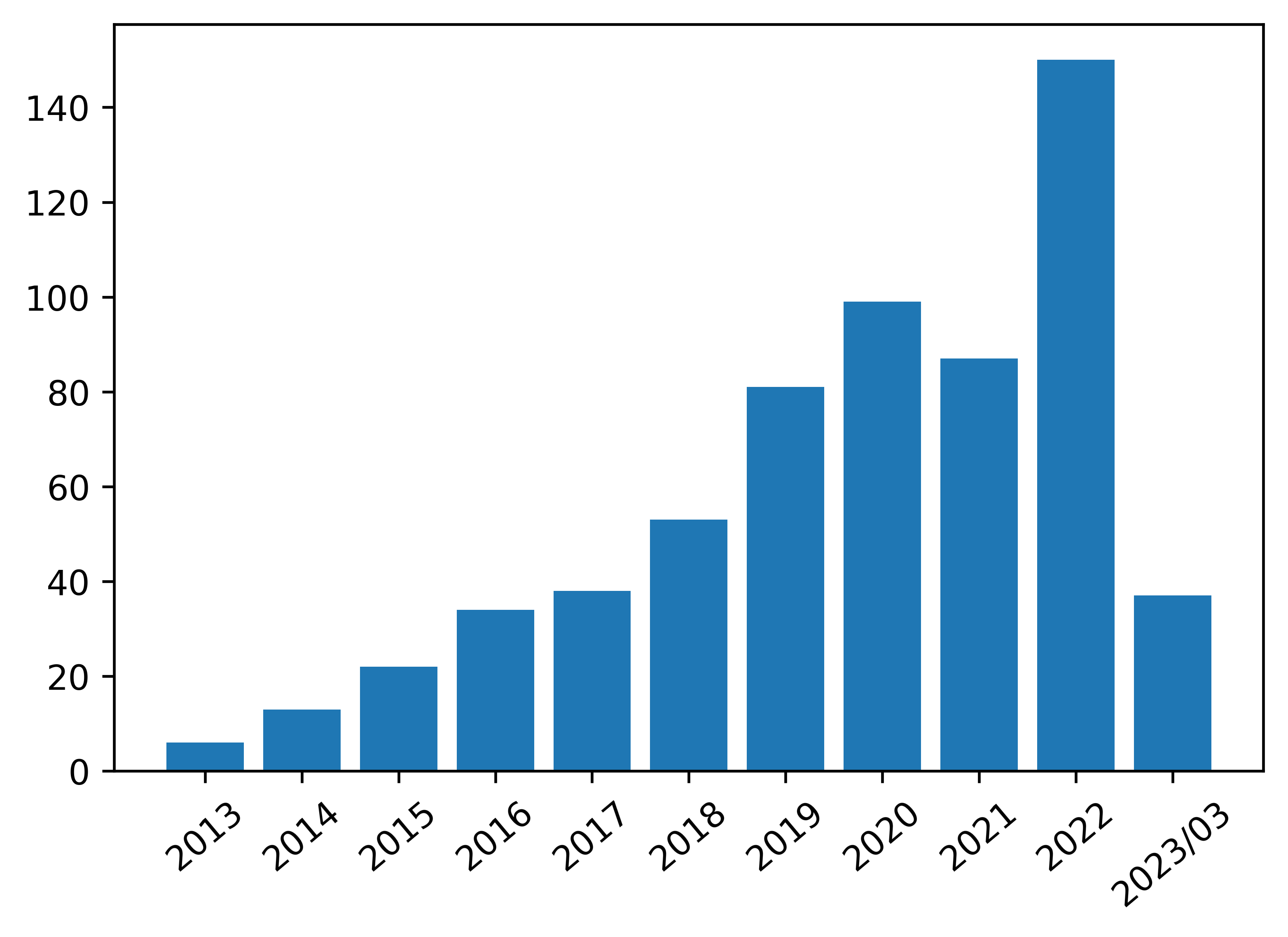}
\caption{Annual distribution of papers}
\label{fig:years}
\end{center}
\end{figure}

Figure \ref{fig:years} shows the number of publications between 2013 and the first quarter of 2023. As seen, there is a strong boost in the number of EES publications in recent last years. This number move from 6 publications in 2013 to reach almost 100 in 2020 and around 90 publications in 2021. it is easy to note the increase in interest in this research field. The average added number of publications was 8 between 2013 and 2017 and around 20.33 in the period 2017-2020. Indeed, this growing interest could be related to many reasons, like energy crises and price fluctuation which put "energy-efficient in manufacturing scheduling" as a trend in the last few years. 

\subsection{Geographical distribution of publications}
\begin{figure}[htbp] 
  \centering
  \includegraphics[width=1.0\textwidth]{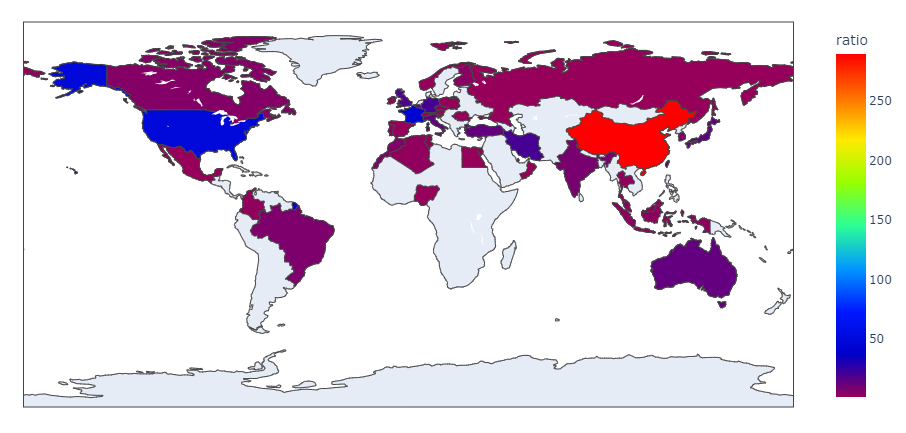}
  \caption{Geographic distribution of publications}
  \label{fig:map2}
\end{figure}

Figure \ref{fig:map2}  presents the geographical distribution of publications, As seen, EES has attracted the interest of many countries with uneven ratios. It is possible to notice four clusters based on the number of publications per country, 
As shown, the first cluster in red represents China which has the biggest part with more than 250 publications in the last 10 years. The second cluster that appears in blue combines the United States and France with around 50 publications each. Then in the third cluster in purple for countries that have around 20 publications such as Germany, the United Kingdom, Iran, and Australia. Finally, many other countries appear in violet like Egypt, Ireland, etc. to represent the fourth class in our focus.

\subsection{Publication distribution among journals and conferences}
The considered papers based on manufacturing scheduling applied to energy efficiency have been published in several journals, books, and conferences with production and operational research orientations. 76 \% of the reviewed documents were published in journals, 23 \% were published in proceedings, and only 1 \% was published as a  book chapter. 
The reviewed publications are related to various topics on scheduling and energy efficiency. Table \ref{Table_journal}, presents the top 15 journals that deal with our topic. Journal of Cleaner Production and International Journal of Production Research are in the foreground respectively with 51 and 27 papers.
Engineering and manufacturing are the main common points between the top journals.

\begin{table}[htbp]
\begin{tabular}{ | p{11cm} |  p{1.5cm} |}
\hline
Journal	&	Number of papers	\\
\hline
Journal of Cleaner Production	&	51	\\
International Journal of Production Research	&	27	\\
Computers and Industrial Engineering	&	21	\\
IEEE Access	&	18	\\
Swarm and Evolutionary Computation	&	13	\\
Sustainability (Switzerland)	&	11	\\
Expert Systems with Applications	&	11	\\
Procedia CIRP	&	10	\\
International Journal of Advanced Manufacturing Technology	&	9	\\
IEEE Transactions on Automation Science and Engineering	&	9	\\
Mathematical Problems in Engineering 	&	7	\\
IFAC-PapersOnLine	&	6	\\
Lecture Notes in Computer Science	&	6	\\
Procedia Manufacturing	&	6	\\
Applied Soft Computing	&	6	\\
\hline 
\end{tabular}
\caption{Distribution of papers published by the journal}
\label{Table_journal}
\end{table}

\section{Literature classification}
\label{sec:classification}
This section describes the main possible classification of reviewed documents
\subsection{Manufacturing setting}
\begin{figure}[hbtp] 
\centering
\includegraphics[width=0.8\textwidth]{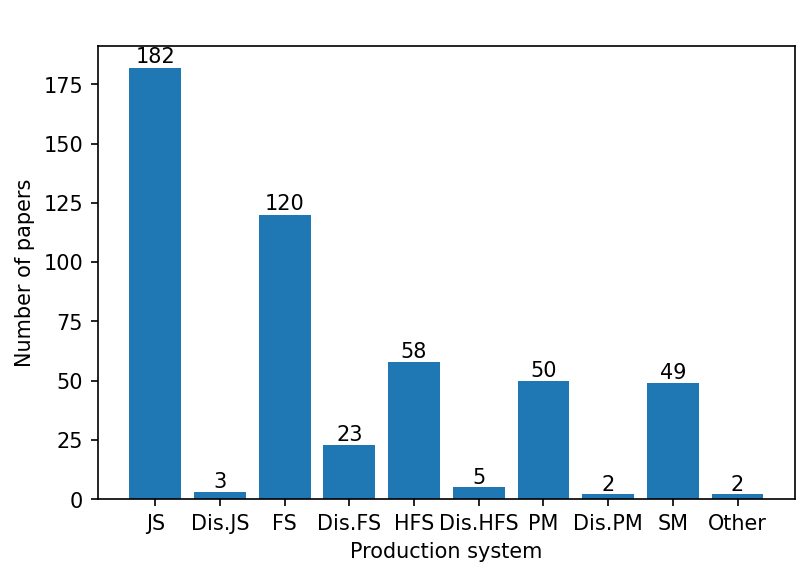}
\caption{\label{fig:prod-sys}Production system distribution}
\end{figure}

The arrangement of machines and the related environmental factors are able to define the category of the production system in order to specify the corresponding model and solving strategies.
As shown in Figure \ref{fig:prod-sys}, 9 classes of production systems are investigated in the literature. Job-shop and flexible job-shop problem accounts for 182 publications as the most considered in the literature such as \cite{Duan2022}, \cite{Salido2017}, \cite{Guo2020}. This kind of production system has attracted the interest of academics and researchers In fact, it is found in all kinds of real-world scenarios including textile, \cite{Ramos2020}, metal industry \cite{coca2019sustainable}. The second investigated production system is the flowshop and flexible flowshop with 120 papers in the last decade. HFS, PM, and SM are also studied with close ratios. Distributed production systems are also considered in energy efficiency, Dis FS is the most tackled with 23 papers while Dis. PM is the least studied with only 2 papers.

\subsection{Theory and practice }
\begin{figure}[hbtp] 
\centering
\includegraphics[width=0.8\textwidth]{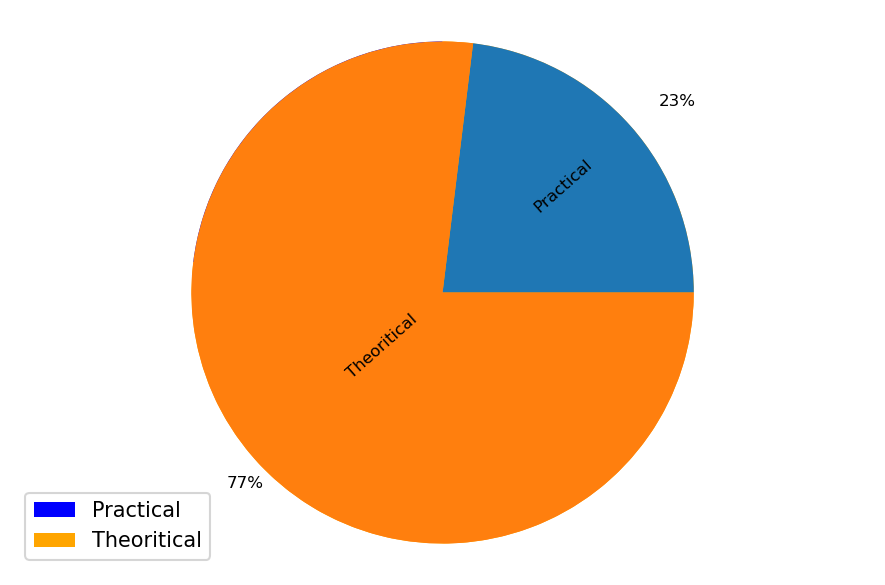}
\caption{\label{fig:thori} Theoretical/practical}
\end{figure}

\begin{figure}[hbtp] 
\centering
\includegraphics[width=1.0\textwidth]{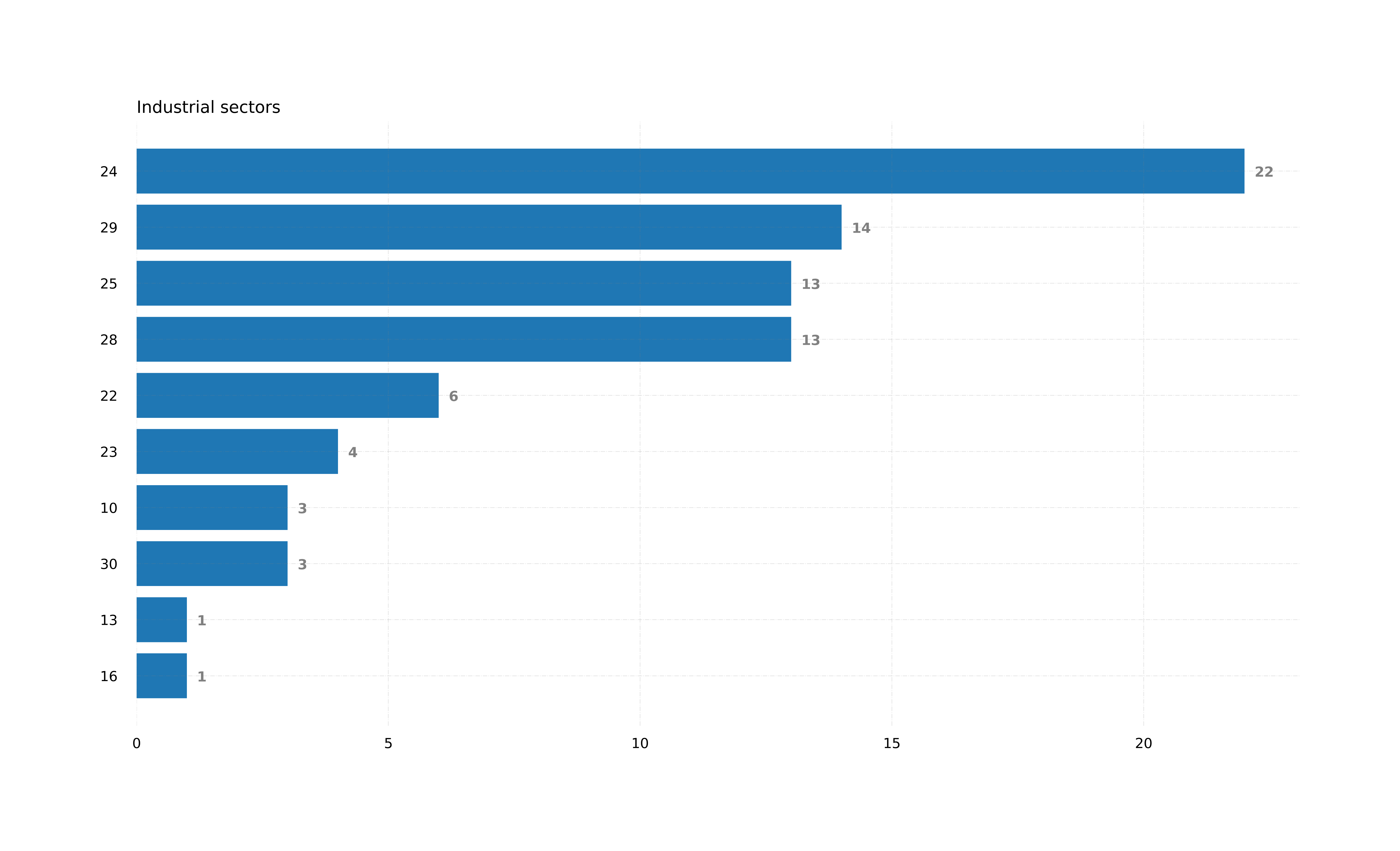}
\caption{\label{fig:sectors} Industrial sectors}
\end{figure}

In literature and according to the type of data we can distinguish between theoretical papers and case studies. As shown in Figure \ref{fig:thori}, 77\% are theoretical papers with randomly generated benchmarks (\cite{Kizilay2019}, \cite{Tan2020}, \cite{xu2014energy}, \cite{Ferretti2020}). In other hand, there is only 23\% of the reviewed papers present practical studies (\cite{Gong2017}, \cite{Liu2020}, \cite{Lu2021}). 
More than 80 papers investigated the EES in practice. All of them are classified in accordance with "International Standard Industrial Classification of All Economic Activities (ISIC 2008)". Figure \ref{fig:sectors}
present the 10 industrial sectors when the basic metal industry (24) is in the foreground and this is very intuitive as it is among the most energy-consuming manufacturing.
In the second level, three sectors are equally considered in the literature, 
Motor vehicles, trailers, and semi-trailers (29),
fabricated metal products, except machinery and equipment (25),
and Manufacture of machinery and equipment n.e.c. (28).
Finally, 6 sectors are less studied in the literature,
Manufacture of rubber and plastics products (22),
Manufacture of other non-metallic mineral products (23)
Manufacture of food products (10),
Manufacture of other transport equipment (30),
Manufacture of textiles (13),
Manufacture of wood and of products of wood and cork, except furniture; manufacture of articles of straw and plaiting materials (16). 
The number of manufacturing sectors in which researchers and academics applied EES is not much compared to the number of illustrated sectors in "International Standard Industrial Classification of All Economic Activities (ISIC 2008)".

\subsection{Number of Objectives}
\begin{figure}[hbtp] 
\centering
\includegraphics[width=0.8\textwidth]{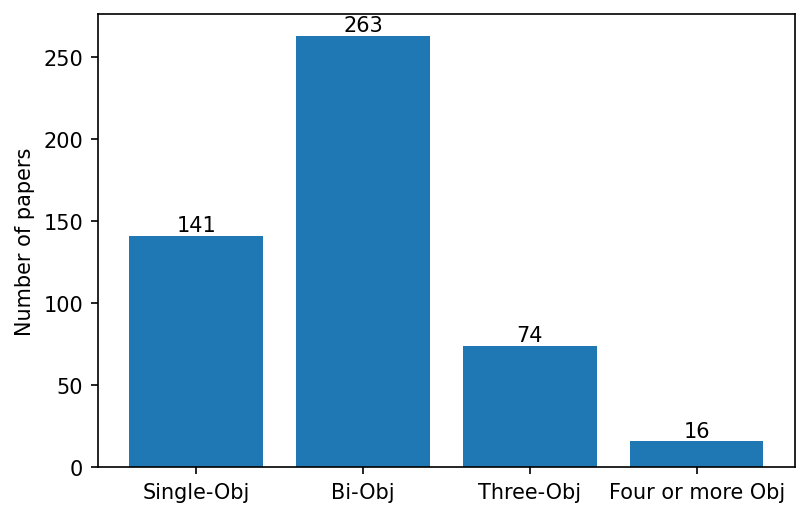}
\caption{\label{fig:number_obj} Distribution of objective function modeling}
\end{figure}

Numerous objective functions have been considered in energy-efficient production scheduling, however, we identified four main categories based on the number of objectives.
These four categories are single objective (\cite{Pang2014}, \cite{Peng2021}, \cite{Nanthapodej2021}), be-objective 
( \cite{wang2020energy}, \cite{Liu2017}, \cite{Wang2016} ), three objectives (\cite{He2022a}, \cite{Faccio2019}, \cite{Wei2021}), and finally, four or more objectives (\cite{Gong2020}). The objectives could be productivity objectives or energy-efficiency objectives or a mixture of them for the three last classes. 
Figure \ref{fig:number_obj} describes the distribution of objective function numbers in the literature of EES.

\begin{figure}[hbtp] 
\centering
\includegraphics[width=0.8\textwidth]{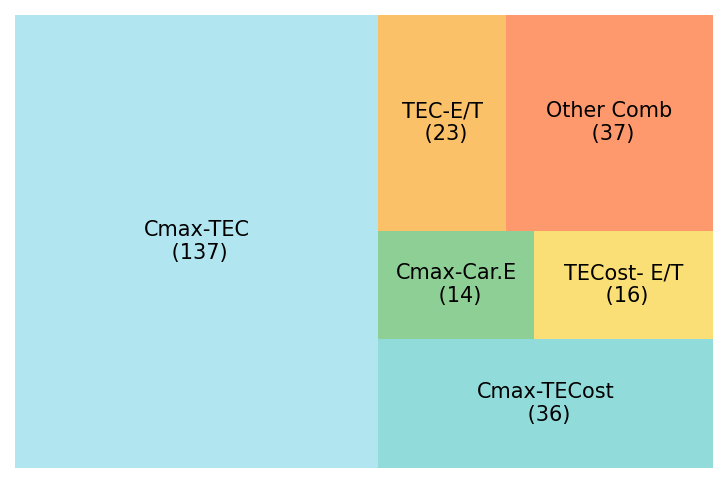}
\caption{\label{fig:bi_obj} Bi-objective distribution}
\end{figure}

\begin{figure}[hbtp] 
\centering
\includegraphics[width=0.8\textwidth]{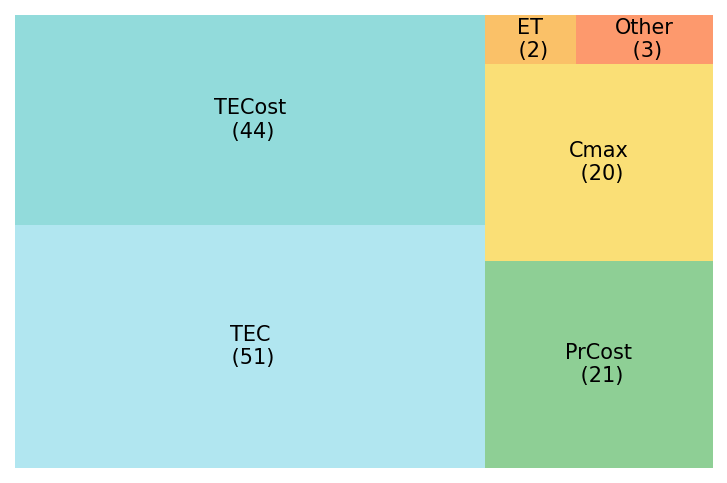}
\caption{\label{fig:single_obj} Single objective distribution}
\end{figure}

The bi-objective category is the most studied with 263 papers. As shown in Figure \ref{fig:bi_obj}, Several combinations of objectives are combined and most of them consider makespan as a productivity measure while TEC and TECost are tackled as energy efficiency measures.
In the second level, the single objective category is studied 141 times in the literature of EES, therefore, most of them are considered to reduce TEC or TECost and CE while other works used productivity measures (Cmax, E/T, etc.) when energy constraint such as peak power is defined in advance. Figure \ref{fig:single_obj} describe better the distribution of single objective functions.
Three objectives problems are tackled 74 times when authors tried to solve more than one productivity or energy efficiency objective at the same time. The last category represents EES problems that combine more than three objectives. Since there are many considered issues in the two last categories, we do not attempt to enumerate all combinations.

\subsection{Energy-efficiency aspects}
Energy-efficient scheduling can be performed by taking into account different aspects involving decisions about energy.

One of the most common ways to minimize energy costs is to schedule jobs at a certain time of the day when the energy price is lower \cite{Kurniawan2021, Aghelinejad2018}. Another strategy is to apply tiered prices in one period \cite{Peng2021, Xu2022}. Time of use and tiered prices are considered respectively in 98 and 2 papers. This is valid in many countries, where those strategies are defined in order to encourage customers to use energy when there is less demand.

Another way to achieve energy efficiency is to use machines with different speeds \cite{Schulz2020, Wei2022, Jiang2019w}. Each speed allows processing a job with a certain processing time and a certain energy consumption, such that the shorter the processing time, the higher the energy consumption, and vice versa. This aspect is investigated 142 in the reviewed papers.

Energy efficiency could be reached by turning off/on, this strategy is considered in 45 papers and it is able to save energy when there is no processing \cite{Lee2017, Liu2019, Gong2021}. However, shutting down and restarting machines could cause energy consumption but it can lead to reducing machine age. In this regard, the energy consumed during turning on/off is considered in 33 papers. 

We mention a typical energy-related constraint, that is the power peak constraint \cite{Masmoudi2019}. This is usually considered when the single objective function is a productivity measure and the power peak constraint sets a limit on the amount of energy used in a certain time instant.

Although Renewable energy is able to provide affordable electricity, this aspect was considered in energy efficient scheduling only in 16 papers in the literature \cite{Subramanyam2020, Wu2020}. 
In those works, wind and solar energies are used in addition to nonrenewable energy to operate manufacturing. In those studies, mathematical models and machine learning are used to deal with uncertainty and predict the amount of energy available during the processing periods.

A very general case consists of optimizing different operating and non-operating modes for the machines which can take several forms.
In operating mode, processing energy consumption is considered in most papers and neglected in a few papers when non-processing energy 
is the objective \cite{ gonzalez2015multi, Peng2018}. 
On the other hand, the non-operating mode of machines can take many forms according to the production system category and constraints.

\begin{figure}[hbtp]
  \centering
  \includegraphics[width=\textwidth]{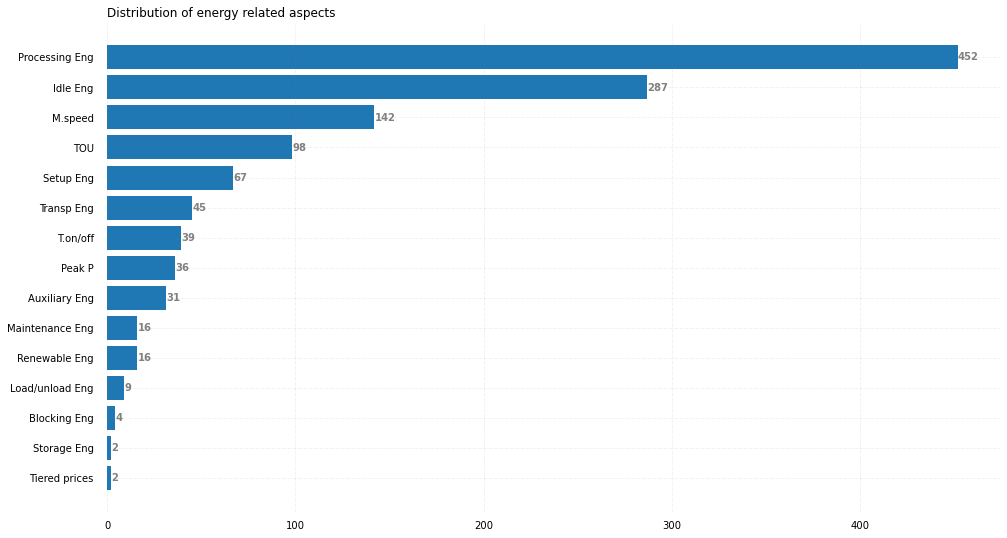}
  \caption{Energy-related aspects}
  \label{fig:aspects}
\end{figure}

Idle energy was the most studied aspect when 287 papers tried to improve its energy consumption \cite{Ambrogio2020, Kawaguchi2019}.
hence, the amount of energy consumed during standby of machines could be important if there is a problem with scheduling or machine allocation.  

Setup times and transportation energy consumption are taken into account 67 and 45 times respectively, those two non-operating tasks are also requiring energy in the production system, hence reducing setup and transportation effort leads to improving energy efficiency \cite{Lu2017, Li2021a}.

A cost in terms of energy is associated with each of these modes and the time spent by a machine in a certain mode is taken into consideration in the computation of the total energy consumption.

Operating and calibrating production lines require common energy. This includes light and air conditioners, and it is not related directly to operating but it is necessary. This issue is considered in 31 papers and called auxiliary energy \cite{Qu2022, Dai2019}.

Figure \ref{fig:aspects} shows the most studied aspects since 2013. The inclusion of operating modes in the problem definition has attracted the most attention from researchers: strategies like the reduction of the machines' idle time are a focal point to avoid energy waste.

\subsection{Other constraints}
\begin{figure}[hbtp] 
\centering
\includegraphics[width=1.0\textwidth]{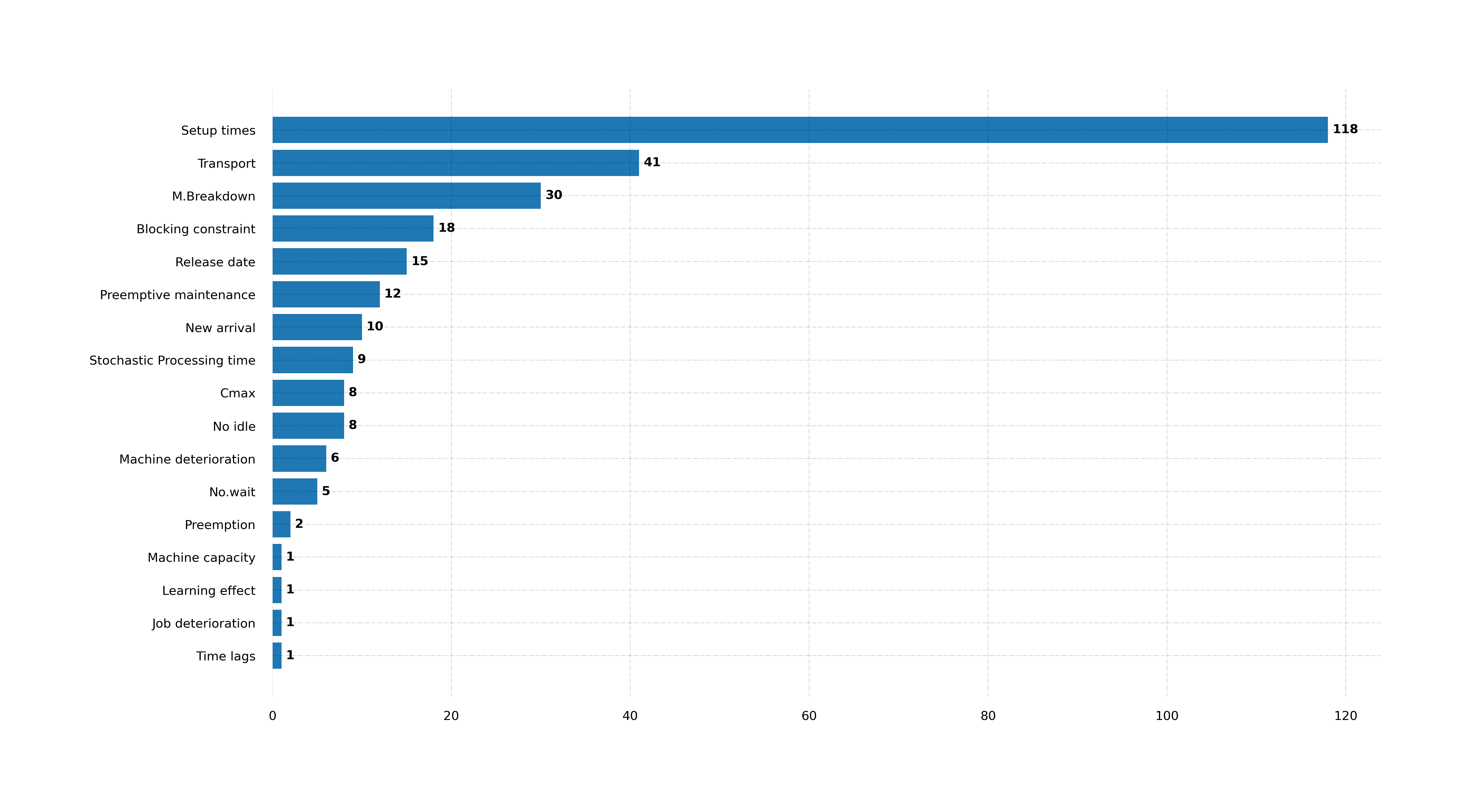}
\caption{\label{fig:constraints} Distribution of scheduling constraints}
\end{figure}

As for the machine environments, we can replicate any scheduling constraints in a study about energy-efficient production scheduling. Here follows a non-exhaustive list of constraints considered since 2013. First, a common constraint for flow shop environment is to impose that all the jobs have to be processed by all the machines in the same order (permutation constraint). 
Figure \ref{fig:constraints} show the distribution of main constraints in EES literature.
Another important constraint is to consider sequence-dependent setup time for the machines \cite{Ramezanian2019}. The authors in \cite{Wu2020} consider the no-wait constraint, meaning that the jobs' waiting time is equal to zero which is essential for process industries. Machine downtime is taken into account in \cite{Park2022}. A variety of process-dependent operational constraints have been also studied on multiple occasions.

\section{Problem formulation and resolution}
\label{sec:formulation}
\subsection{Model formulation}
\begin{figure}[hbtp] 
\centering
\includegraphics[width=0.8\textwidth]{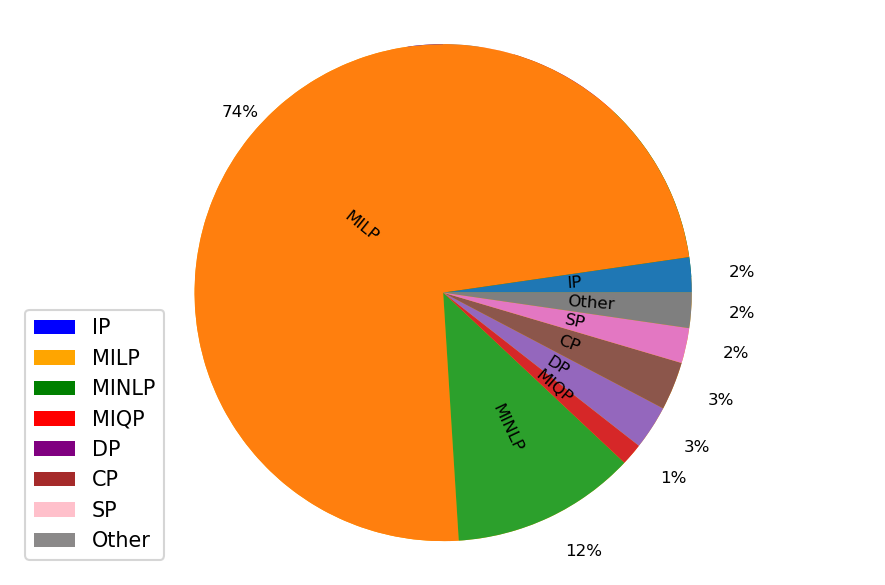}
\caption{\label{fig:math_model} Problem formulation}
\end{figure}

Agnostic to the manufacturing system setting, many papers provide a formal model for the problem they study. Figure \ref{fig:math_model} describes the use of the different modeling approaches. Mixed integer linear programming (MILP) is the most used in 74\% of papers that consider mathematical models \cite{fathollahi2021sustainable, Cheng2022} and at a lower rate, Mixed integer non-linear programming MINLP is used in 12\% \cite{Lee2017, Sin2020}. Many other approaches are used with small ratios such as Mixed integer quadratic programming MIQP \cite{Chen2022}, Constraint programming CP \cite{Ham2021}, Dynamic programming DP \cite{Heydar2022} etc. 

\subsection{Objective function formations}
\begin{figure}[hbtp] 
\centering
\includegraphics[width=0.8\textwidth]{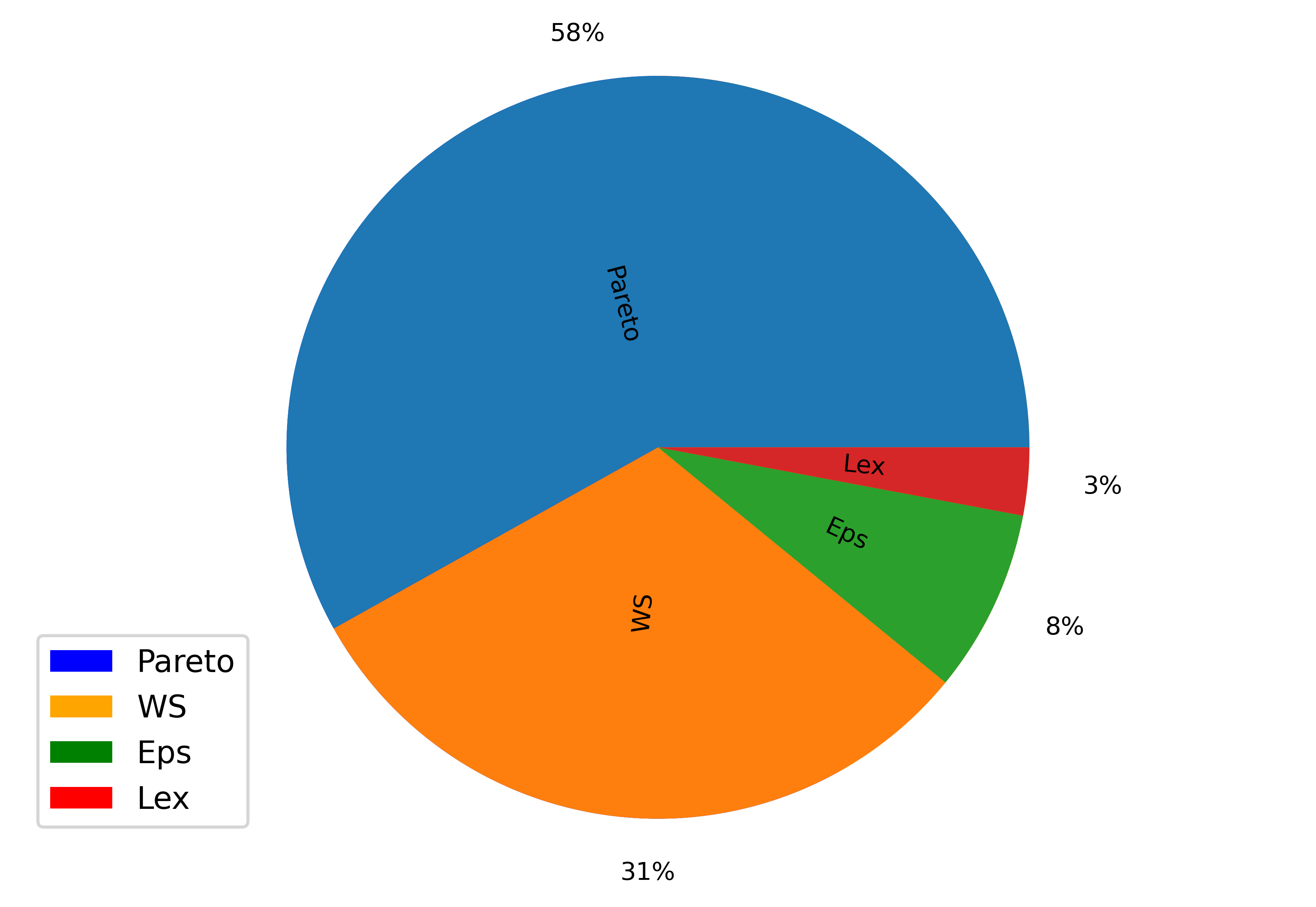}
\caption{\label{fig:obj_for} Objective function formations}
\end{figure}

Many objectives are taken into consideration in EES literature, 
Those objectives are investigated in order to optimize productivity (makespan and earliness tardiness) or energy efficiency (TEC, TECost, etc.). In most papers, the improvement of one objective could affect the other and hence cause conflicts. This conflict is resolved by finding the optimal balance between objectives. Beyond the question of the number of considered objectives, they could be formed in many ways. Figure \ref{fig:obj_for} show the considered multiobjective formation considered in the last 10 years according to \cite{t2006multicriteria}.
Pareto dominance is the most used tool with a ratio of 58\%, this method does not require prior preference knowledge about the criterion from decision-makers. applying this method results in a  Pareto front containing all the best solutions for multiobjective problems \cite{Dai2019}.
In other cases, the decision-maker could Determine the importance of each criterion and assign a weight to each of them. this technique is the weighted sum (WS) is considered at 31\% of the multiobjective reviewed papers \cite{Lv2022}. 
Finally, epsilon constraint (eps) and lexicographical (lex) are used respectively in 8\% and 3\% of the reviewed papers.

\subsection{Resolution Methodology}
As for all combinatorial problems, energy-efficient scheduling problems can be addressed by both exact and approximation approaches. The selection of the methods is related to many factors such as the instance size, the number of objectives, and the problem structure. However, exact approaches are not performing well for practical problems, especially in larger instances. In this work, we adopted a simple classification with three known classes, exact, heuristics, and metaheuristics.

\subsubsection{Exact methods}
Exact methods are able to give the optimal solution for a given problem. Branch and bound ($B \& B$) is the most advisable technique for solving EES problems optimally. Exact approaches are rarely used in the literature of EES given the computational complexity.
\cite{Assia2022} proposed a ($B \& B$) approach to solving unrelated parallel machines with Cmax and TEC objectives, the problem is modeled as MILP and solved using a commercial solver. In \cite{Cui2020}, the ($B \& B$) is used to tackle the single\-machine problem under the time\-of\-use electricity tariff in order to minimize Cmax and power costs. A commercial solver is used to compute the Pareto solution. \cite{assia2020genetic} solved the single machine under the constraints of system robustness and stability in order to minimize TEC using ($B \& B$), the problem was modeled with MILP and solved using a solver. Dynamic programming is also used for exact resolution in \cite{Cui2021}, this approach is proposed to minimize TECost under the TOU tariff in the FS problem
Finally, the exact decomposition approach has been used in a few research works to handle complicated constraints and solve energy-efficient scheduling problems to optimal \cite{Liu2017b}.
Despite the relative success, exact algorithms are still inefficient with large instances.
\subsubsection{Heuristics}
Dispatching rules are widely used in scheduling literature and are the simplest type of heuristics. Those methods provide some simple rules that are able to rank jobs and assign them to machines.
\cite{Guo2021} used the shortest processing time rule (SPT) and adaptive genetic algorithm to tackle the single-machine problem to minimize the total waiting time aiming to achieve energy savings.
In \cite{Zandi2020}, authors investigated the problem of the parallel machine with Cmax and TEC aiming to obtain the
exact Pareto frontier using a heuristic algorithm. 
In \cite{Li2016a}, the unrelated parallel machines problem is considered where total tardiness and TEC are considered as objective functions. Three kinds of constructive heuristics are used for the resolution.
\cite{Wang2018a} proposed a new constructive heuristic to study the
 two-stage hybrid flowshop scheduling with Cmax criteria and energy-saving considerations. 
In \cite{Ho2021}, six constructive heuristics are used to tackle the flowshop with TECost criteria.

\subsubsection{Meta-heuristics}
\begin{figure}[hbtp] 
\centering
\includegraphics[width=0.8\textwidth]{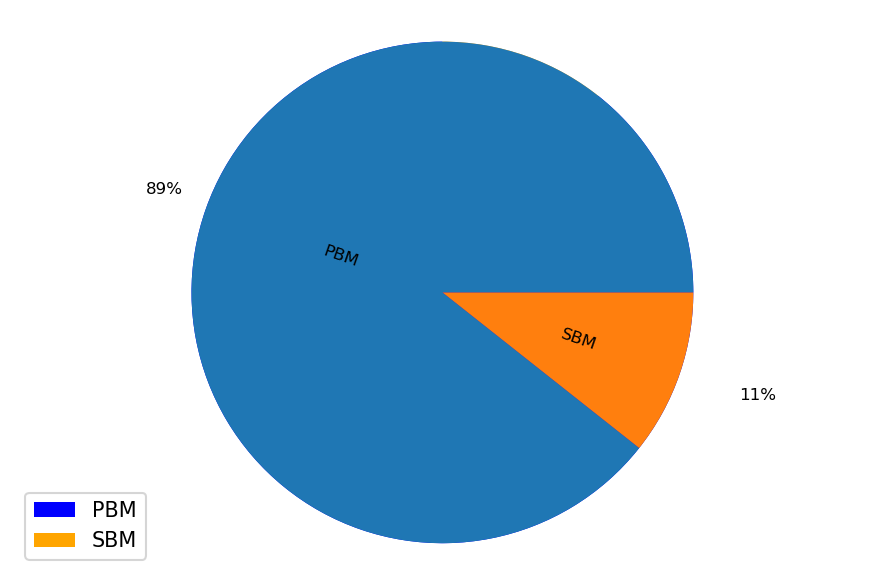}
\caption{\label{metaheuristics} Type of metaheuristics}
\end{figure}

In the last few decades, researchers and academics have developed many strategies to improve the efficiency of simple heuristics. 
The idea is to start with an initial solution, add randomness to the search, and increase replications in order to obtain better solutions. This kind of method is called metaheuristics. We classified these methods according to the type of initial solution. The use of these two kind of metaheuristics is given in Figure \ref{metaheuristics}
\begin{itemize}
    \item Single solution-based metaheuristics
\end{itemize}
This category starts with a single initial solution and progresses in the search space by describing a trajectory to improve the objective function value. 
In \cite{Nanthapodej2021}, authors proposed a variable (VNS) neighborhood search to minimize TEC in parallel machines scheduling problems with Job Priority and Control completion time.
The HFS scheduling problem is tackled in \cite{Keller2015} for TEC minimization using improved simulated annealing (SA) 
In \cite{Cheng2021}  authors proposed an iterated greedy (IG) to minimize the weighted sum of the makespan and total flowtime value in a no-idle flowshop with energy savings. As shown in Figure \ref{metaheuristics}, SBM represents only 11\% of the used metaheuristics, most of them are used to study single objective EES problems.

\begin{itemize}
    \item Population-based metaheuristics
\end{itemize}
Those procedures start with more than one candidate solution called the population of solution. In this category, the most used methods are related to Swarm Intelligence and Evolutionary Computation algorithms. 
Evolutionary procedures including genetic algorithms (GA) and memetic (MA) algorithms are investigated in many cases for the resolution of different production systems. 
\cite{ning2021low} an improved GA is considered to solve the flexible JS with carbon emission minimization.
\cite{yaurima2018hybrid} a bi-objective HFS is tackled with a non-dominated Sorting Genetic Algorithm II (NSGA-II).
\cite{Lu2021} a MA is proposed to solve the Dist-FS with  Cmax, negative social impact (NSI), and TEC minimization.
In \cite{jiang2019improved}, the authors proposed a multi-objective evolutionary algorithm to study the FS problem with setup times.
Swarm Intelligence metaheuristics have been used frequently
in EES literature. Particle swarm intelligence (PSO) is used in \cite{Kawaguchi2018, Abedi2017, Liu2019}. Ant colony optimization (ACO) has been used in \cite{Zheng2020, Qiao2022}. 
As shown in Figure \ref{metaheuristics}, 89\% of the used procedures are PBM.

\section{Analysis of the literature }\label{sec:interaction}
Many objectives are considered in EES literature, some parts of them are investigated to improve energy issues while the other part is to deal with productivity. To achieve those two kinds of objectives, we have to consider many energy aspects that reflect the general environmental factors and the nature of objective functions. This section is to gather insights about energy-related aspects and the addressed objectives as they engage with objectives categories.

\subsection{Association between objective functions and energy aspects}
\begin{figure}[hbtp] 
\centering
\includegraphics[width=0.8\textwidth]{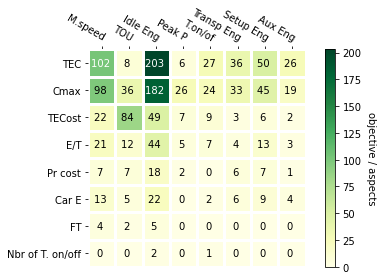}
\caption{\label{fig:heatmap_obj_aspect} Interaction between Objective function and energy aspects}
\end{figure}

Figure \ref{fig:heatmap_obj_aspect} presents the heatmap of the main objectives function and the most used energy aspects in the last 10 years. A strong picture is given by the heatmap to visualize the frequency of each energy aspect for each objective. We can clearly notice that TEC is the objective function, Idle energy is considered in 203 cases and the speed level of machines is taken into account in 102 papers in a lesser way setup energy and transportation energy are taken into consideration respectively in 50 and 36 times. Also as shown, there is a high correlation between TECost and TOU and in a lesser way with idle energy. For the Cmax objective, many aspects are present with uneven ratios, the most frequent of them are machine speed, idle energy, and setup energy. Peak power is considered when Cmax is studied. Most other objectives are correlated more with idle energy.

\subsection{Association between objective functions and number of objectives}
\begin{figure}[hbtp] 
\centering
\includegraphics[width=0.8\textwidth]{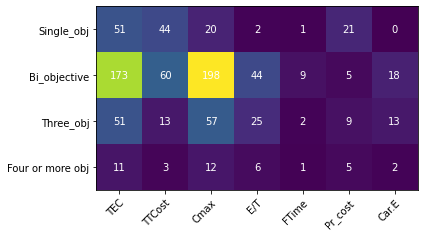}
\caption{\label{fig:heatmap_obj_cat} Distribution of objectives function per category}
\end{figure}

To understand the interest of decision-makers in EES, a heatmap was drawn to show the presence of the main objectives among their four categories. In the plot \ref{fig:heatmap_obj_cat}, notice the predominance of TEC improvement as an energy efficiency measure in single and double objectives while the TECost is less addressed especially when talking about three or more objectives.
Moving to productivity measures, Cmax is the most studied objective for all categories, this criterion 
draws its importance from the ability to maximize the use of available resources which reduces idle times and therefore idle energy, also Cmax could reduce auxiliary energy by finishing as soon as possible the production process. Many other criteria are addressed with uneven ratios such as E/T and FT.

\section{Future work and conclusions}
\label{sec:futurework}

The growing heed in sustainability in the last few years gave rise to increased attention for integrating sustainability aspects into manufacturing. moreover, energy has a substantial impact on manufacturing companies and their economic performance, therefore, it is important to improve energy efficiency.
It is intuitive that the improving efficiency of machinery is a well-known method to achieve more efficient and sustainable production processes up to a certain extent. However, a system-level approach for scheduling manufacturing jobs has a significant potential to improve energy efficiency (\cite{carbonfoot}). The biggest benefit of a system-level approach is the quick win without requiring a big investment. As a result, there has been a significant increase in the number of EES publications taking into consideration many kinds of production systems with diverse objectives, several constraints, and energy aspects. Consequently, the multiplicity of terms and notation in this field makes their classification challenging. In this paper, we introduced a systematic literature review taking into account all scopes and aims in the last 10 years.
The first insight from our work is providing a guideline to researchers with the most promising area along with their statistical and geospatial analysis.  

%In the second step, this study presented a framework to classify 506 papers according to many criteria such as production system, energy-related aspects, number of objectives, industrial sectors, and resolution methodology. The framework is able to help academics and researchers who are interested in EES.  

One of the important insights we conclude in this paper is the limited number of industrial applications. There is a big opportunity to apply presented theoretical work in practice. From another practical perspective, we observed that the impact of emerging energy supply chain issues is not covered in detail. There is certainly a big potential for researchers in this domain. 

From a technical point of view, we see a clear distinction in the literature in terms of exploiting model formulations. Considerably a few studies considered exploiting mathematical and artificial intelligence (particularly constraint programming) methods for EES. Indeed, the declarative approach of these methods can be integrated with the fast search capability of heuristic/meta-heuristic methods to achieve better quality EES solutions in shorter times. 

Another technical opportunity in EES is integrating solution methods discussed in this paper with simulation techniques in which more realistic and high-fidelity solutions can be obtained (\cite{simheuristic}). 

With given energy and supply crisis along with the inflationist global economy, authors expect that EES will stay one of the most active research areas in combinatorial optimization literature in the near future.

\bibliographystyle{unsrtnat}
\bibliography{references}  %%% Uncomment this line and comment out the ``thebibliography'' section below to use the external .bib file (using bibtex) .

%%% Uncomment this section and comment out the \bibliography{references} line above to use inline references.
% \begin{thebibliography}{1}

% 	\bibitem{kour2014real}
% 	George Kour and Raid Saabne.
% 	\newblock Real-time segmentation of on-line handwritten arabic script.
% 	\newblock In {\em Frontiers in Handwriting Recognition (ICFHR), 2014 14th
% 			International Conference on}, pages 417--422. IEEE, 2014.

% 	\bibitem{kour2014fast}
% 	George Kour and Raid Saabne.
% 	\newblock Fast classification of handwritten on-line arabic characters.
% 	\newblock In {\em Soft Computing and Pattern Recognition (SoCPaR), 2014 6th
% 			International Conference of}, pages 312--318. IEEE, 2014.

% 	\bibitem{hadash2018estimate}
% 	Guy Hadash, Einat Kermany, Boaz Carmeli, Ofer Lavi, George Kour, and Alon
% 	Jacovi.
% 	\newblock Estimate and replace: A novel approach to integrating deep neural
% 	networks with existing applications.
% 	\newblock {\em arXiv preprint arXiv:1804.09028}, 2018.

% \end{thebibliography}

\end{document}